\documentclass[11pt, reqno]{amsart}
\usepackage{amssymb}
\usepackage{amscd}
\usepackage{verbatim,ifthen}
\usepackage{color}
\usepackage{latexsym}

\usepackage{mathrsfs}
\usepackage{wrapfig}

\usepackage{color}

\usepackage{bbm}

\usepackage{array}
\newcolumntype{M}[1]{>{\centering\arraybackslash}m{#1}}

\usepackage[colorlinks, linkcolor=black, citecolor=blue, linktocpage,urlcolor=blue]{hyperref}
\addtocontents{toc}{\setcounter{tocdepth}{1}}

\usepackage{amsmath}

\numberwithin{equation}{section}

\usepackage{eucal}
\usepackage{calc}  
\usepackage{enumitem} 
\usepackage{tensor}
\usepackage{graphicx,wrapfig,lipsum}
\usepackage{etoolbox}

\makeatletter
\patchcmd{\@mn@margintest}{\@tempswafalse}{\@tempswatrue}{}{}
\patchcmd{\@mn@margintest}{\@tempswafalse}{\@tempswatrue}{}{}
\reversemarginpar 
\makeatother
\usepackage{scrextend}

\makeatletter
\DeclareRobustCommand\widecheck[1]{{\mathpalette\@widecheck{#1}}}
\def\@widecheck#1#2{%
    \setbox\z@\hbox{\m@th$#1#2$}%
    \setbox\tw@\hbox{\m@th$#1%
       \widehat{%
          \vrule\@width\z@\@height\ht\z@
          \vrule\@height\z@\@width\wd\z@}$}%
    \dp\tw@-\ht\z@
    \@tempdima\ht\z@ \advance\@tempdima2\ht\tw@ \divide\@tempdima\thr@@
    \setbox\tw@\hbox{%
       \raise\@tempdima\hbox{\scalebox{1}[-1]{\lower\@tempdima\box
\tw@}}}%
    {\ooalign{\box\tw@ \cr \box\z@}}}
\makeatother
\title{A General Schwarz Lemma For Hermitian Manifolds}
\author{Kyle Broder}
\address{The University of Queensland,  St. Lucia,  QLD 4067, Australia}
\email{k.broder@uq.edu.au}
\thanks{The first named author was supported by funding from the Australian Government through the Australian Research Council's Discovery Projects funding scheme (project DP220102530).  The second named author was supported by an Australian Government Research Training Program (RTP) Scholarship.}
\author{James Stanfield}
\address{The University of Queensland,  St. Lucia,  QLD 4067, Australia}
\email{james.stanfield@uq.net.au} 
\usepackage{soul}

\begin{document}

\maketitle

\begin{abstract}
The Schwarz lemma for holomorphic maps between Hermitian manifolds is improved.  New curvature constraints on the source and target manifolds are introduced and shown to be weaker than the Ricci and real bisectional curvature, respectively.  The novel target curvature condition is intrinsic to the Hermitian structure and is controlled by the holomorphic sectional curvature if the metric is pluriclosed.  This leads to significant improvements on the Wu--Yau theorem.  Further, it is shown that the Schwarz lemma is largely invariant under a change of Hermitian connection, and the precise geometric quantity that varies with the connection is determined.  This enables us to establish the Schwarz lemma for the Gauduchon connections.
\end{abstract}

\section{Introduction and Main Results}
A compact complex manifold $X$ is \textit{Kobayashi hyperbolic} if every holomorphic map $\mathbf{C} \to X$ is constant. Such manifolds exhibit a very rich structure and have maintained a central position in complex geometry since their discovery (see, e.g., \cite{KobayashiBook}). Despite the unrelenting attention they have received, a complete understanding of the canonical bundle $K_X$ of such manifolds has remained out of reach. A (folklore generalization of) a long-standing conjecture made by Kobayashi predicts that a compact Kobayashi hyperbolic manifold is projective with ample canonical bundle. In particular, by the Aubin--Yau theorem \cite{Aubin,Yau1976},  the Kobayashi conjecture predicts that all compact Kobayashi hyperbolic manifolds admit a unique K\"ahler--Einstein metric with negative scalar curvature. An important source of Kobayashi hyperbolic manifolds comes from compact complex manifolds with a Hermitian metric of negative \textit{holomorphic sectional curvature} $\text{HSC}_{\omega} <0$, where $$\text{HSC}_{\omega}(\zeta) \ := \ \frac{1}{| \zeta |_{\omega}^4} \sum_{\alpha, \beta, \gamma, \delta} R_{\alpha \overline{\beta} \gamma \overline{\delta}} \zeta^{\alpha} \overline{\zeta}^{\beta} \zeta^{\gamma} \overline{\zeta}^{\delta}, \hspace*{1cm} \zeta \in T^{1,0}X.$$ The (Tosatti--Yang extension \cite{TosattiYang} of the) Wu--Yau theorem \cite{
WuYau1} states that a compact K\"ahler manifold $(X,\omega)$ with a K\"ahler metric satisfying $\text{HSC}_{\omega}<0$ is projective with ample canonical bundle, verifying the Kobayashi conjecture and a conjecture of S.-T. Yau (see, e.g., \cite{HeierLuWong}) for this class of manifolds. To address the canonical bundle of general compact Kobayashi hyperbolic manifolds, however, it is essential to consider non-K\"ahler Hermitian metrics. Even on compact K\"ahler manifolds, if the metric with $\text{HSC}_{\omega}<0$ is not K\"ahler, then the ampleness of $K_X$ remains wide open. In this setting, Yang--Zheng \cite{YangZhengRBC} introduced the real bisectional curvature $\text{RBC}_{\omega}$, whose negativity is sufficient to deduce that $K_X$ is ample if $X$ is compact K\"ahler. Until now, there have been no non-trivial extensions of the Wu--Yau theorem to compact complex manifolds with $\text{HSC}_{\omega}<0$. 

The first main theorem exhibits the most general form of the Wu--Yau theorem on compact K\"ahler manifolds, requiring that the negatively curved metric is \textit{pluriclosed} in the sense that $\sqrt{-1} \partial \bar{\partial} \omega=0$.

\subsection*{Theorem 1.1}\label{thm:main}
Let $X$ be a compact K\"ahler manifold with a pluriclosed metric satisfying $\operatorname{HSC}_\omega < 0$. Then $X$ is projective with ample canonical bundle. In particular, $X$ admits a unique K\"ahler--Einstein metric of negative scalar curvature.\\

\nameref{thm:main} is a consequence of a more general statement (see \nameref{MainWuYau}) which holds for general Hermitian metrics, under a stronger curvature constraint that intermediates between $\text{HSC}_{\omega}<0$ and $\text{RBC}_{\omega}<0$. The crux of the proof of the Wu--Yau theorem is the Schwarz lemma \cite{Lu,Royden,YauSchwarzLemma, YangZhengRBC},  the primary technique for estimating how two metrics relate geometrically from assumptions on their curvature. In more detail, the Schwarz lemma is the result of applying the maximum principle to a formula for the Laplacian of the energy density $\Delta | \partial f |^2$ of a holomorphic map $f :(X,\omega_g) \to (Y, \omega_h)$.  The computation of $\Delta_{\omega_g} | \partial f |^2$ for a general holomorphic map was first carried out by Lu \cite{Lu}, where it was shown that \begin{eqnarray}\label{LaplacianEnergyDensity}
\Delta_{\omega_g} | \partial f |^2 &=& | \nabla \partial f |^2 + \left( \text{Ric}_{\omega_g}^{(2)} \right)^{\sharp}  \otimes \omega_h(\partial f, \overline{\partial f}) - \omega_g^{\sharp} \otimes \omega_g^{\sharp} \otimes  R^h(\partial f, \overline{\partial f}, \partial f, \overline{\partial f}),
\end{eqnarray} where $\omega_g^{\sharp}$ is the induced metric  on the cotangent bundle $\Omega_X^1$ (c.f., \cite[$\S 3$]{BroderSBC}). In \eqref{LaplacianEnergyDensity}, the \textit{second Chern Ricci curvature} is defined ${\text{Ric}_{\omega_g}^{(2)}= \sqrt{-1} \sum_{k,\ell} \text{Ric}_{k \bar{\ell}}^{(2)} e^k \wedge \bar{e}^{\ell}}$, whose components in a local frame are given by $\text{Ric}_{k \bar{\ell}}^{(2)} = g^{i \bar{j}} R_{i \bar{j} k \bar{\ell}}$. 

Despite the vast number of applications and attention the Schwarz lemma has received,  there have been only a limited number of general (i.e., not specific to a given context) improvements in obtaining estimates from \eqref{LaplacianEnergyDensity}. The two most notable achievements are due to Yau \cite{YauSchwarzLemma} and Royden \cite{Royden}. Yau's breakthrough \cite{YauSchwarzLemma} was the application of his generalized maximum principle \cite{YauMaximumPrinciple, Omori} which enabled the Schwarz lemma to be applied to a significantly larger class of holomorphic maps. On the other hand, Royden \cite{Royden} showed that if the target metric $\omega_h$ is K\"ahler,  then an upper bound on the holomorphic sectional curvature of $\omega_h$ is sufficient to control the \textit{target curvature term} $\omega_g^{\sharp} \otimes \omega_g^{\sharp} \otimes  R^h(\partial f, \overline{\partial f}, \partial f, \overline{\partial f})$.

Yang--Zheng \cite{YangZhengRBC} carried out a systematic investigation of the target curvature term,  in the absence of any K\"ahler assumption.  They showed that it can be expressed as  \begin{eqnarray}\label{RBCDefine}
\text{RBC}_{\omega}(\xi) \ : = \ \frac{1}{| \xi |_{\omega}^2} \sum_{\alpha, \beta, \gamma, \delta}  R_{\alpha \overline{\beta} \gamma \overline{\delta}} \xi^{\alpha \overline{\beta}} \xi^{\gamma \overline{\delta}},
\end{eqnarray} for non-negative Hermitian $(1,1)$--tensors $\xi$,  and referred to \eqref{RBCDefine} as the \textit{real bisectional curvature}. Yang--Zheng showed that the real bisectional curvature dominates the holomorphic sectional curvature in general, and is comparable to it (in the sense that they always have the same sign) when the metric is K\"ahler.

\nameref{thm:main} follows from the second main theorem of the present article, which provides a significant general improvement on the Schwarz lemma for holomorphic maps between Hermitian manifolds. Let us define, for $\tau \in \mathbf{R}_{\geq 0}$, the \textit{tempered real bisectional curvature} \begin{eqnarray}\label{TemperedRBCDefinitionEquation}
\text{RBC}_{\omega_h}^{\tau}(\xi) &=& \frac{1}{| \xi |_{\omega}^2} \sum_{\alpha, \beta, \gamma, \delta} \left( R_{\alpha \bar{\beta} \gamma \bar{\delta}} - \frac{1-\tau}{4} \sum_{\rho, \sigma} T_{\alpha \gamma}^{\rho} \overline{T_{\beta \delta}^{\sigma}} h_{\rho \overline{\sigma}} \right) \xi^{\alpha \bar{\beta}} \xi^{\gamma \bar{\delta}},
\end{eqnarray} where $\xi$ is a non-negative Hermitian $(1,1)$--tensor, and $T$ is the torsion of the Chern connection. In the special case that $\tau=1$,  we recover the \textit{real bisectional curvature} $\text{RBC}_{\omega_h}$, but for $\tau<1$, negative tempered real bisectional curvature $\text{RBC}_{\omega_h}^{\tau}<0$ is decidedly weaker than $\text{RBC}_{\omega}<0$ (see \nameref{RBCExample}). If the Hermitian metric is pluriclosed, then for $\tau=0$, the tempered real bisectional curvature is comparable to the holomorphic sectional curvature (see \nameref{PluriclosedLemma}).

For $\tau \in \mathbf{R}_{>0} \cup \{ \infty \}$, we define the \textit{tempered Ricci curvature} \begin{eqnarray}\label{TemperedRicciFormula}
\text{Ric}_{\omega_g}^{\tau}  &:= & \text{Ric}_{\omega_g}^{(2)} + \frac{1-1/\tau}{4}\mathcal{Q}_{\omega_g}^{\circ},
\end{eqnarray} where $\text{Ric}_{\omega_g}^{(2)}$ is the second Chern Ricci curvature and $\mathcal{Q}_{\omega_g}^{\circ} = \sqrt{-1} \sum_{k,\ell} \mathcal{Q}_{k\bar{\ell}}^{\circ} e^k \wedge \bar{e}^{\ell}$, where $\mathcal{Q}_{k\bar{\ell}}^{\circ} = \sum_{p,q,s,r} T_{pr}^i \overline{T_{qs}^j} g_{k\bar{j}} g_{i\bar{\ell}}$ (c.f., \cite[Notation 3.3]{BroderStanfield}). From \eqref{TemperedRicciFormula}, if $\tau=1$, the tempered Ricci curvature recovers the second Chern Ricci curvature. The tempered Ricci curvature $\text{Ric}_{\omega_g}^{\tau}$ serves to atone for the crimes that we commit on the target manifold  by introducing the tempered real bisectional curvature $\text{RBC}_{\omega_h}^{\tau}$ in the main general Schwarz lemma:

\subsection*{Theorem 1.2}\label{MainTheoremTemperedRBC}
Let $f : (X, \omega_g) \to (Y, \omega_h)$ be a holomorphic map of rank $r$ between Hermitian manifolds.  Suppose there are constants $\tau \in (0,\infty)$, $C_1 \in \mathbf{R}$ and $C_2, \kappa_0 \geq 0$ such that $$\text{Ric}_{\omega_g}^\tau \ \geq \ -C_1 \omega_g + C_2 f^{\ast} \omega_h, \hspace{1cm}  \text{RBC}_{\omega_h}^\tau \leq - \kappa_0 \leq 0.$$ Then \begin{eqnarray}\label{SchwarzInequality1}
\Delta_{\omega_g} | \partial f |^2 & \geq & -C_1 | \partial f |^2 + \left( \frac{\kappa_0}{r} + \frac{C_2}{n} \right) | \partial f |^4. 
\end{eqnarray} In particular, if in addition, $X$ is compact and $\kappa_0 n + r C_2 >0$, then \begin{eqnarray}\label{SchwarzLemma0}
| \partial f |^2 & \leq & \frac{C_1 r n}{\kappa_0 n + r C_2}.
\end{eqnarray}

\hfill

\noindent The above theorem can be extended to non-compact manifolds, so long as $\omega_g$ is complete with (Riemannian) Ricci curvature bounded from below and bounded torsion $(1,0)$--form. The proof of \nameref{MainTheoremTemperedRBC} makes use of a careful analysis of the Hessian term in \eqref{LaplacianEnergyDensity}. We use precisely the skew-symmetric part of $| \nabla \partial f|^2$ to temper the source and target curvature terms using the torsion. 

The last main theorem we discuss here indicates that \nameref{MainTheoremTemperedRBC} is close to optimal.  There has long been interest in producing a Schwarz lemma that involves the curvature (and torsion) of general Hermitian connections, not exclusively the Chern connection ${}^c \nabla$.  The most notable class of connections being the \textit{Gauduchon connections} \cite{GauduchonHermitianConnections} \begin{eqnarray}\label{GaudFormula}
{}^t \nabla &: =& t {}^c \nabla + (1-t) {}^{\ell } \nabla,
\end{eqnarray}  where $t \in \mathbf{R}$ and ${}^{\ell} \nabla$ is the \textit{Lichnerowicz connection} (the Hermitian part of the Levi-Civita connection). The geometry of such connections has received considerable attention (see, e.g., \cite{GauduchonHermitianConnections, BroderStanfield,LafuenteStanfield,StanfieldThesis} and the references therein).

The authors \cite{BroderStanfield} discovered that the Gauduchon holomorphic sectional curvature ${}^t \text{HSC}_{\omega}$ exhibits a monotonicity phenomenon ${}^t \text{HSC}_{\omega} \leq {}^c \text{HSC}_{\omega}$, for all $t \in \mathbf{R}$.  It is, therefore, natural to suspect that improvements in the Schwarz lemma could be achieved by varying the connection. The following result shows, however, that (for a more general class of so-called `suitably skew-symmetric Hermitian connections', that include the Gauduchon connections, see \nameref{Def52}), the Schwarz lemma estimate obtained from discarding the Hessian term, is independent of the connection:

\subsection*{Theorem 1.3} \label{thm:OtherConnections}
Let $f : (X, \omega_g) \to (Y, \omega_h)$ be a holomorphic map between Hermitian manifolds. Let $\nabla$ and $\widetilde{\nabla}$ be suitably skew-symmetric Hermitian connections on $T^{1,0}X$ and $T^{1,0}Y$, and write ${}^{\mathcal{T}} \nabla$ for the connection on $\mathcal{T} : = \Lambda_X^{1,0} \otimes f^{\ast} T^{1,0}Y$ induced by $\nabla$ and $\widetilde{\nabla}$. Then \begin{eqnarray*}
\Delta_{\omega_g} | \partial f |^2 &=& \left | \text{Sym}\left( {}^{\mathcal{T}} \nabla^{1,0} \partial f \right) \right|^2 + |\partial f \circ {}^c T - {}^c \widetilde{T}(\partial f \cdot,\partial f \cdot)|^2\\&&\qquad+ \left( {}^c \text{Ric}_{\omega_g}^{(2)} \right)^{\sharp} \otimes \omega_h(\partial f, \overline{\partial f}) - \omega_g^{\sharp} \otimes \omega_g^{\sharp}  \otimes {}^c\widetilde{R}(\partial f, \overline{\partial f}, \partial f, \overline{\partial f}),
\end{eqnarray*} where ${}^c T$ and ${}^c \widetilde{T}$ denotes the Chern torsion of $\omega_g$ and $\omega_h$, respectively, and ${}^c \widetilde{R}$ denotes the Chern curvature of $\omega_h$.

\hfill

As a consequence of \nameref{thm:OtherConnections}, we establish the Schwarz lemma for the Gauduchon connections (see \nameref{Prop44}), a result that has long been sought. Given the formidable nature of the calculation and formula, let us mention the important special case of the Strominger--Bismut connection (corresponding to the Gauduchon connection $t=-1$):

\subsection*{Theorem 1.4}
Let $f : (X, \omega_g) \to (Y, \omega_h)$ be a rank $r$ holomorphic map from a compact Hermitian manifold to a pluriclosed manifold with $\text{HSC}_{\omega_h} \leq - \kappa_0 <0$.  Suppose that there are constants $C_1 \in \mathbf{R}$, $C_2 \geq 0$ such that \begin{eqnarray*}
2 \left( {}^b \text{Ric}_{\omega_g}^{(1)} + \sqrt{-1}\text{Re} {}^b T_{kr}^{\ell} \overline{{}^b T_{ir}^i} e^k \wedge \overline{e^{\ell}} \right)  - {}^b \text{Ric}_{\omega_g}^{(2)} + {}^b \text{Ric}_{\omega_g}^{(3)} + {}^b \text{Ric}_{\omega_g}^{(4)} & \geq & - C_1 \omega_g + C_2 f^{\ast} \omega_h.
\end{eqnarray*} Then there is a constant $\kappa_1>0$ such that \begin{eqnarray*}
\Delta_{\omega_g} | \partial f |^2 & \geq & - \frac{C_1}{3} | \partial f |^2 + \left( \frac{C_2}{3n} + \frac{\kappa_1}{r} \right) | \partial f |^4.
\end{eqnarray*} In particular,  if $C_2 r + 3n \kappa_1 > 0$,  then \begin{eqnarray*}
| \partial f |^2 & \leq & \frac{C_1 nr }{C_2 r + 3n \kappa_1}.
\end{eqnarray*}

\subsection*{Acknowledgements}
The first named author would like to thank Ben Andrews, Artem Pulemotov, and Gang Tian for their valuable support and encouragement. He would also like to thank Franc Forstneri\v{c}, Herv\'{e} Gaussier, Ramiro Lafuente, Finnur L\'arusson, Peter Petersen, Miles Simon, and Stefano Trapani for their interest in this work, and valuable discussions. The second named author would like to thank Ramiro Lafuente for engaging discussions and advice.

\section{The Tempered Curvatures}\label{sec:TC}
Let $X$ be a complex manifold with underlying complex structure $J$. We identify $\omega_g(\cdot, \cdot):= g(J \cdot, \cdot)$ with the underlying Hermitian metric $g$. We denote by $\nabla = {}^c \nabla$ the Chern connection, the unique Hermitian connection such that $\nabla^{0,1} = \bar{\partial}$. If $\{ e_i \}$ is a local frame for the tangent bundle $T^{1,0}X$, we write $T_{ij}^k e_k : = T(e_i, e_j)$ for the components of the Chern torsion $$T(u,v) \ := \ \nabla_u v - \nabla_v u - [u,v].$$ We write $R_{i \bar{j}k \bar{\ell}} := R(e_i, \bar{e}_j, e_k, \bar{e}_{\ell})$ for the components of the Chern curvature tensor $$R(u,\bar{v}, w, \bar{z}) \ := \ g \left(\nabla_u \nabla_{\bar{v}} w - \nabla_{\bar{v}} \nabla_u w - \nabla_{[u,\bar{v}]} w, \bar{z} \right).$$

We start by showing that for $\tau=0$, the tempered real bisectional curvature $\text{RBC}_{\omega}^{\tau}$, defined in \eqref{TemperedRBCDefinitionEquation}, is intrinsic to the Hermitian structure.

\subsection*{Lemma 2.1}\label{CoordinateLemma}
Let $(X, \omega_g)$ be a Hermitian manifold.  For any point $p \in X$, there are local holomorphic coordinates centered at $p$ such that for all $1 \leq i,j,k,\ell \leq n$, $$g_{k \overline{\ell}} \ = \ \delta_{k\overline{\ell}}, \hspace{1cm} \partial_i g_{k \overline{\ell}} \ = \ \frac{1}{2} T_{ik}^{\ell}, \hspace*{1cm} \partial_i  \partial_{\overline{j}} g_{k \overline{\ell}} \ = \ - R_{i \overline{j} k \overline{\ell}} + \frac{1}{4} T_{ik}^p \overline{  T_{j\ell}^q} g_{p \overline{q}}.$$ 

\noindent In particular, the tempered real bisectional curvature is the second-order correction term for the Hermitian metric in the geodesic normal coordinates of the Chern connection.  \begin{proof}
Recall that (see, e.g.,  \cite[Lemma 2.9]{StreetsTian}) in the geodesic normal coordinates coming from the Chern connection,  the components $g_{k \overline{\ell}}$ of the metric satisfy $g_{k\overline{\ell}}(p) = \delta_{k\overline{\ell}}(p)$ and $\partial_i g_{k \overline{\ell}}(p) = - \partial_k g_{i \overline{\ell}}(p)$. In any local frame,  the torsion and curvature of the Chern connection are given by \begin{eqnarray*}
T_{ij}^k &=& g^{k \overline{\ell}} \left( \partial_i g_{j \overline{\ell}} - \partial_j g_{i \overline{\ell}} \right), \hspace*{1cm} R_{i \overline{j} k \overline{\ell}} \ = \ - \partial_i \partial_{\overline{j}} g_{k \overline{\ell}} + g^{p \overline{q}} \partial_i g_{k \overline{q}} \partial_{\overline{j}} g_{p \overline{\ell}}.
\end{eqnarray*} Hence,  in the geodesic normal coordinates at the point $p \in X$,  the torsion reads $T_{ij}^k(p) = 2 \partial_i g_{j \overline{\ell}}$.  Similarly, the tempered curvature tensor $R_{i \overline{j} k \overline{\ell}} - \frac{1}{4} T_{ik}^p \overline{T_{j\ell}^q}$ is given in these coordinates by \begin{eqnarray*}
R_{i \overline{j} k \overline{\ell}} - \frac{1}{4} T_{ik}^p \overline{T_{j\ell}^q} &=&  - \partial_i \partial_{\overline{j}} g_{k \overline{\ell}} + \partial_i g_{k \overline{p}} \partial_{\overline{j}} g_{p \overline{\ell}} - \frac{1}{4} T_{ik}^p \overline{T_{j\ell}^p} \\
&=&  - \partial_i \partial_{\overline{j}} g_{k \overline{\ell}} + \frac{1}{4} T_{ik}^p \overline{T_{j\ell}^p} - \frac{1}{4} T_{ik}^p \overline{T_{j\ell}^p} \ = \  - \partial_i \partial_{\overline{j}} g_{k \overline{\ell}},
\end{eqnarray*}
as required.
\end{proof}

The following example illustrates that the negativity of the tempered real bisectional curvature is, at least locally, strictly weaker than the negativity of the real bisectional curvature.

\subsection*{Example 2.2}\label{RBCExample}
Let $A \in \mathbf{C}^n \otimes \Lambda^2(\mathbf{C}^n)^{\star}$. For $\varepsilon>0$, define a Hermitian metric on the ball $\mathbf{B}_r(0) \subset \mathbf{C}^n$, for $r>0$ sufficiently small, by the formula \begin{eqnarray*}
g_{k \overline{\ell}} & : = & \delta_k^{\ell} + A_{ik}^{\ell} z^i + \overline{A_{ik}^{\ell} z^i} + \frac{1}{2} A_{ik}^p \overline{A_{j\ell}^p} z^i \overline{z}^j + \varepsilon z^{\ell} \overline{z}^k.
\end{eqnarray*} 

At $z=0$, we have $T_{ij}^k = 2 A_{ij}^k$ and $R_{i \overline{j} k \overline{\ell}} = \frac{1}{2} A_{ik}^p \overline{A_{j\ell}^p} - \varepsilon \delta_k^j \delta_i^{\ell}$.  Moreover, $R_{i \overline{j} k \overline{\ell}} - \frac{1}{4} T_{ik}^p \overline{T_{j\ell}^q} = - \frac{1}{2} A_{ik}^p \overline{A_{j\ell}^p} - \varepsilon \delta_k^j \delta_i^{\ell}$.  For $r$ sufficiently small, the metric has negative tempered real bisectional curvature $ \text{RBC}_{\omega_g}^{0} <0$, but for $A\neq 0$ and $\varepsilon \ll | A |^2$,  the real bisectional curvature will not be negative.  \\

An essential component of the proof of \nameref{thm:main} is showing that if the metric $\omega$ is pluriclosed, then $\text{HSC}_{\omega}<0$ implies that $\text{RBC}_{\omega}^{\tau}<0$ are equivalent for $\tau=0$. This is achieved in the following lemma and illustrates one of the essential advantages of the improved Schwarz lemma that is not exhibited by previous forms of the Schwarz lemma.

\subsection*{Lemma 2.3} \label{PluriclosedLemma}
Let $(X, \omega)$ be a pluriclosed manifold.  Then $\text{HSC}_{\omega}$ and $\text{RBC}_{\omega}^{\tau}$ for $\tau=0$ are comparable in the sense that they always have the same sign. In particular, $\text{HSC}_{\omega} \leq 0$ (respectively, $\text{HSC}_{\omega}<0$) implies $\text{RBC}_{\omega}^0 \leq 0$ (respectively, $\text{RBC}_{\omega}^0 <0$). \begin{proof}
Let $\omega$ be a pluriclosed metric. Then from the first Bianchi identity,  the pluriclosed condition $\sqrt{-1} \partial \bar{\partial} \omega =0$ is equivalent to \begin{eqnarray}\label{PluriclosedSymmetry}
R_{\alpha \bar{\beta} \gamma \bar{\delta}} - R_{\gamma \bar{\beta} \alpha \bar{\delta}} - R_{\alpha \bar{\delta} \gamma \bar{\beta}} + R_{\gamma \bar{\delta} \alpha \bar{\beta}} &=& T_{\alpha \gamma}^{\rho} \overline{T_{\beta \delta}^{\sigma}} h_{\rho \overline{\sigma}},
\end{eqnarray} see, e.g.,  \cite{YeDerivativeEst} for the details.  From \eqref{PluriclosedSymmetry}, we observe that \begin{eqnarray*}
\frac{1}{2} \sum_{\alpha, \beta, \gamma, \delta} \left( R_{\alpha \bar{\beta} \gamma \bar{\delta}} - R_{\alpha \bar{\delta} \gamma \bar{\beta}} \right) \xi^{\alpha \bar{\beta}} \xi^{\gamma \bar{\delta}} &=& \frac{1}{4} \sum_{\alpha, \beta, \gamma, \delta, \rho, \sigma} T_{\alpha \gamma}^{\rho} \overline{T_{\beta \delta}^{\sigma}} h_{\rho \overline{\sigma}} \xi^{\alpha \bar{\beta}} \xi^{\gamma \bar{\delta}}.
\end{eqnarray*} The holomorphic sectional curvature is comparable (in the sense that it always has the same sign) to the \textit{altered holomorphic sectional curvature} \begin{eqnarray*}
\widetilde{\text{HSC}}_{\omega}(\xi) &=& \frac{1}{| \xi |_{\omega}^2} \sum_{\alpha, \beta, \gamma, \delta} \left( R_{\alpha \bar{\beta} \gamma \bar{\delta}} + R_{\alpha \bar{\delta} \gamma \bar{\beta}} \right) \xi^{\alpha \bar{\beta}} \xi^{\gamma \bar{\delta}},
\end{eqnarray*} which was formally defined in \cite{BroderTangAltered,BroderStanfield}, but appeared earlier in the literature implicitly (c.f., \cite{YangZhengRBC}).  We compute the tempered real bisectional curvature for $\tau=0$ when the metric is pluriclosed: \begin{eqnarray*}
\text{RBC}_{\omega}^{\tau}(\xi) &=& \frac{1}{| \xi |_{\omega}^2} \sum_{\alpha, \beta, \gamma, \delta} \left( R_{\alpha \bar{\beta} \gamma \bar{\delta}} - \frac{1}{4} \sum_{\rho, \sigma} T_{\alpha \gamma}^{\rho} \overline{T_{\beta \delta}^{\sigma}} h_{\rho \overline{\sigma}} \right) \xi^{\alpha \bar{\beta}} \xi^{\gamma \bar{\delta}} \\
&=& \frac{1}{| \xi |_{\omega}^2} \sum_{\alpha, \beta, \gamma, \delta} \left( R_{\alpha \bar{\beta} \gamma \bar{\delta}} - \frac{1}{2} \left( R_{\alpha \bar{\beta} \gamma \bar{\delta}} - R_{\alpha \bar{\delta} \gamma \bar{\beta}}  \right) \right) \xi^{\alpha \bar{\beta}} \xi^{\gamma \bar{\delta}} \ = \ \frac{1}{2} \widetilde{\text{HSC}}_{\omega}(\xi).
\end{eqnarray*}
\end{proof}

We conclude this section with a discussion of Hermitian metrics with pointwise constant tempered real bisectional curvature. The result will further illustrate the difference between the real bisectional curvature and the tempered real bisectional curvature. Yang--Zheng \cite[Theorem 2.9]{YangZhengRBC} showed that if the real bisectional curvature is pointwise constant $\text{RBC}_{\omega} = \kappa_0$, then $\kappa_0 \leq 0$. Moreover, if $\kappa_0=0$, then all Chern Ricci curvatures vanish and the metric is \textit{balanced} in the sense that the torsion $(1,0)$--form $\eta$ defined by \begin{eqnarray}\label{torsion1form}
\partial \omega^{n-1} &=& -\eta \wedge \omega^{n-1},
\end{eqnarray} vanishes. For more results of this type, we invite the reader to see \cite{BroderTangAltered} and the references therein. For the tempered real bisectional curvature, we have the following significantly more general result.

\subsection*{Proposition 2.4}
Let $(X, \omega)$ be a compact Hermitian manifold with $\dim_{\mathbf{C}} X=n$. Fix $\tau \in \mathbf{R}_{\geq 0}$. If the tempered real bisectional curvature is pointwise constant $ \text{RBC}_{\omega}^\tau = \kappa_0$ for some $\kappa_0 \in \mathbf{R}$, then \begin{eqnarray*}
2\kappa_0n(n-1)\text{Vol}_{\omega}(X) &=& (1-\tau)\int_X |T|^2 \omega^n  - 2\int_X |\eta|^2 \omega^n
\end{eqnarray*} In particular, if $\tau<1$ and $\omega$ is balanced, then $\kappa_0 \leq 0$ and the metric is K\"ahler. If $\tau>1$ and $n \geq 2$, then $\kappa_0 \leq 0$ with $\kappa_0 =0$ if and only if the metric is K\"ahler.

\begin{proof}
Following \cite{YangZhengRBC,BroderTangAltered}, let \( \{ e_i \} \) be a unitary frame. If the tempered real bisectional curvature is pointwise constant \( \text{RBC}_{\omega}^\tau = \kappa_0 \), then 
\begin{eqnarray}\label{ConstantEqns}
    R_{i\overline j k \overline \ell} + R_{k \overline \ell i \overline j} \ = \ 2 \kappa_0 \delta_{i\ell} \delta_{kj} + \frac{1-\tau}{2}\sum_r T_{ik}^r\overline{T^r_{j\ell}},
\end{eqnarray}
for all indices \( 1\leq i,j,k,\ell \leq n \). Let \( \{ e^k \} \) be a unitary coframe of \( (1,0) \)--forms dual to the unitary frame \( \{ e_k \} \). With respect to this local frame, the torsion $(1,0)$--form is given by \( \eta = \sum_{i,j}T_{ij}^i e^j \). From the first Bianchi identity \( T_{ij, \overline{\ell}}^k = R_{j \overline{\ell} i \overline{k}} - R_{i \overline{\ell} j \overline{k}} \), we see that 
\begin{eqnarray*}
\eta_{j,\overline{\ell}} &=& \sum_k (R_{j \overline{\ell} k \overline{k}} - R_{k \overline{\ell} j \overline{k}}),
\end{eqnarray*} 
where the index after the comma indicates covariant differentiation with respect to the Chern connection. From \eqref{ConstantEqns}, we compute 
\begin{eqnarray}
\sum_i \eta_{i, \overline{i}} &=& \sum_{i,k} (R_{i \overline{i} k \overline{k}} - R_{k \overline{i} i \overline{k}}) \nonumber \\
&=& \sum_{i,k} \left(-R_{k\overline k i \overline i} + R_{i \overline k k \overline i} + 2\kappa_0 \delta_k^i \delta^k_i - 2 \kappa_0 \delta^i_i\delta_k^k + (1-\tau)\sum_r |T^r_{ik}|^2\right) \nonumber \\
&=& -\sum_i\eta_{i, \overline i} + 2\kappa_0n(1-n) + (1-\tau)|T|^2. \label{IntegrateProof}
\end{eqnarray} Differentiating \eqref{torsion1form}, we have
\[
\overline{\partial} \partial (\omega^{n-1}) \ = \ \bar{\partial} \eta \wedge \omega^{n-1} + \eta \wedge \bar{\eta} \wedge \omega^{n-1}.
\] 
Hence,  by integrating over $X$ and using \eqref{IntegrateProof}, we have
\[
2 \int_X | \eta |^2 \omega^n \ = \ 2 \int_X \sum_i \eta_{i, \overline{i}} \omega^n \ = \  2\kappa_0n(1-n) \text{vol}_{\omega}(X) + (1-\tau)\int_X |T|^2 \omega^n,
\] where $\text{vol}_{\omega}(X) : = \int_X \omega^n$.
\end{proof}

\section{Proofs of \nameref{thm:main} and \nameref{MainTheoremTemperedRBC}} \label{sec:mainProofs}
In this section, we prove the novel Schwarz lemma estimate \nameref{MainTheoremTemperedRBC} together with the application to the Wu--Yau theorem in \nameref{thm:main}.

\subsection*{Proof of \nameref{MainTheoremTemperedRBC}} 
Let $f : (X, \omega_g) \to (Y, \omega_h)$ be a rank $r$ holomorphic map. We denote by $T, R$ and $\widetilde{T}, \widetilde{R}$ the Chern torsion and curvature of $\omega_g$ and $\omega_h$, respectively. Fix a point $x \in X$ and let $e_i$ and $\widetilde{e}_{\alpha}$ be local frames at $x \in X$ and $f(x)\in Y$, respectively. Write $\partial f(e_i) = f_i^{\alpha} \widetilde{e}_{\alpha}$. A calculation going back to Lu \cite{Lu} gives \begin{eqnarray}\label{eqn:BochnerPartialf}
\Delta_{\omega_g} | \partial f |^2 &=& |  \widehat{\nabla}^{1,0} \partial f |^2 + \text{Ric}_{q \overline p}^{(2)} f_p^{\alpha} \overline{f_q^{\alpha}} - \widetilde{R}^h_{\alpha \overline{\beta} \gamma \overline{\delta}}  f_q^{\alpha} \overline{f_q^{\beta}} f_p^{\gamma} \overline{f_{p}^{\delta}},
\end{eqnarray}
where $\widehat \nabla$ denotes the connection on $\Lambda_X^{1,0} \otimes f^{\ast} T^{1,0}Y$ induced from the Chern connections of $\omega_g$ and $\omega_h$. The symmetric part $\operatorname{Sym}\left( \widehat \nabla^{1,0} \partial f\right)$ of the Hessian $\widehat \nabla^{1,0} \partial f$ in these frames reads
\[
\operatorname{Sym}\left(\widehat \nabla^{1,0} \partial f \right)_{ij}^\alpha := \tfrac{1}{2}\left((\widehat \nabla_i \partial f)_j^\alpha -(\widehat \nabla_j\partial f)_i^\alpha \right).
\]
Since the symmetric and skew-symmetric parts are orthogonal, we may write 
\begin{eqnarray*}
| \widehat{\nabla}^{1,0} \partial f |^2 &=& |  \widehat{\nabla}^{1,0} \partial f - \text{Sym}(\widehat{\nabla}^{1,0} \partial f ) |^2 + | \text{Sym}(\widehat{\nabla}^{1,0} \partial f) |^2 \\
&=& \frac{1}{4} \sum_{k, \ell, \alpha} | (\widehat{\nabla}_k \partial f )_{\ell}^{\alpha} - (\widehat{\nabla}_{\ell} \partial f)_k^{\alpha} |^2 + | \text{Sym}(\widehat{\nabla}^{1,0} \partial f ) |^2.
\end{eqnarray*} Let $\Gamma$ and $\widetilde{\Gamma}$ respectively denote the Christoffel symbols for the Chern connection of $\omega_g$ and $\omega_h$. Then  \begin{eqnarray*}
(\widehat{\nabla}_k \partial f)_{\ell}^{\alpha} - (\widehat{\nabla}_{\ell} \partial f)_k^{\alpha} &=& f_k^{\gamma} \widetilde{\Gamma}_{\gamma \rho}^{\alpha} f_{\ell}^{\rho} - \Gamma_{k\ell}^p f_p^{\alpha} - ( f_{\ell}^{\gamma} \widetilde{\Gamma}_{\gamma \rho}^{\alpha} f_k^{\rho} - \Gamma_{\ell k}^p f_p^{\alpha}) \ = \ f_k^{\gamma} f_{\ell}^{\rho} \widetilde{T}_{\gamma \rho}^{\alpha} - T_{k\ell}^p f_p^{\alpha}.
\end{eqnarray*}
Hence, \begin{eqnarray*}
|  \nabla^{1,0} \partial f |^2 &=& \frac{1}{4} \sum_{k,\ell,\alpha} | ( \nabla_k \partial f)_{\ell}^{\alpha} - ( \nabla_{\ell} \partial f)_k^{\alpha} |^2 + | \text{Sym}(\nabla^{1,0} \partial f ) |^2 \\
&=& \frac{1}{4} \sum_{k,\ell, \alpha} |f_k^{\gamma} f_{\ell}^{\rho} \widetilde{T}_{\gamma \rho}^{\alpha} - T_{k\ell}^p f_p^{\alpha}|^2 + | \text{Sym}( \nabla^{1,0} \partial f ) |^2.
\end{eqnarray*} From Young's inequality with $\tau \in (0,\infty)$, \begin{eqnarray*}
\frac{1}{4} \sum_{k,\ell, \alpha} | \widetilde{T}_{\gamma \rho}^{\alpha} f_k^{\gamma} f_{\ell}^{\rho} - T_{k\ell}^p f_p^{\alpha} |^2 & \geq & \frac{1}{4}(1-\tau)  \sum_{k,\ell,\alpha}| \widetilde{T}_{\gamma \rho}^{\alpha} f_k^{\gamma} f_{\ell}^{\rho} |^2 + \frac{1}{4}(1-\tau^{-1}) \sum_{k,\ell,\alpha}|T_{k\ell}^p f_p^{\alpha} |^2.
\end{eqnarray*}
Substituting this estimate into \eqref{eqn:BochnerPartialf} gives \begin{eqnarray*}
\Delta_{\omega_g} | \partial f |^2 & \geq & | \text{Sym}(\widehat{\nabla}^{1,0} \partial f ) |^2  + \left(  \text{Ric}_{k \overline{\ell}}^{(2)} + \frac{1}{4}(1-\tau^{-1}) T_{p q}^\ell \overline{T_{pq}^k} \right) f_\ell^{\alpha} \overline{f_k^{\alpha}}\\
&& \hspace{2cm} - \left( \widetilde{R}_{\alpha \overline{\beta} \gamma \overline{\delta}} - \frac{1}{4}(1-\tau) \widetilde{T}_{\alpha \gamma}^{\mu} \overline{\widetilde{T}_{\beta \delta}^{\mu}}\right)f_p^{\alpha} \overline{f_p^{\beta}} f_q^{\gamma} \overline{f_{q}^{\delta}}.
\end{eqnarray*} Since $\text{Ric}_{\omega_g}^{\tau} \geq - C_1 \omega_g + C_2 f^{\ast} \omega_h$, by choosing the unitary frames such that $f_i^{\alpha} = \lambda_i \delta_i^{\alpha}$, for $\lambda_i \geq 0$, $1 \leq i \leq n$, we have \begin{eqnarray*}
 \text{Ric}_{k \overline{\ell}}^{\tau} f_\ell^{\alpha} \overline{f_k^{\alpha}} & \geq & - C_1 \sum_i\lambda_i^2 + C_2 \sum_i\lambda_i^4 \ \geq \ - C_1 | \partial f |^2 + \frac{C_2}{n} | \partial f |^4.
\end{eqnarray*} Further, the bound, $\text{RBC}_{\omega_h}^{\tau} \leq - \kappa_0$ implies \begin{eqnarray*}
 \left( \widetilde{R}_{\alpha \overline{\beta} \gamma \overline{\delta}} - \frac{1}{4}(1-\tau) \widetilde{T}_{\alpha \gamma}^{\mu} \overline{\widetilde{T}_{\beta \delta}^{\mu}}\right)f_p^{\alpha} \overline{f_p^{\beta}} f_q^{\gamma} \overline{f_{q}^{\delta}} & \leq & - \frac{\kappa_0}{r} | \partial f |^4.
\end{eqnarray*} Combining these estimates gives \eqref{SchwarzLemma0}. The second claim follows from applying the Omori--Yau maximum principle \cite{Omori, YauMaximumPrinciple} as in \cite{YauSchwarzLemma,YangZhengRBC}. \hfill $\Box$

\hfill

\noindent If the source metric is K\"ahler, we may permit $\tau \searrow 0$. In particular, if $\omega_h$ is Hermitian with $\text{RBC}_{\omega_h}^0 \leq -\kappa_0 < 0$, or pluriclosed with $\text{HSC}_{\omega_h} \leq -\kappa_0 \leq 0$, then may apply \nameref{PluriclosedLemma} and \nameref{MainTheoremTemperedRBC} to obtain:

\subsection*{Corollary 3.1}\label{KahlerNewSL}
Let $f\colon (X,\omega_g) \to (Y,\omega_h)$ be a holomorphic map of rank $r$ from a complete K\"ahler manifold to a Hermitian manifold. Suppose that $\omega_h$ is Hermitian with $\text{RBC}_{\omega_h}^{\tau} \leq - \kappa_0 < 0$, for $\tau \geq 0$, or pluriclosed with $\text{HSC}_{\omega_h} \leq - \kappa_0 \leq 0$. If there are constants $C_1 \in \mathbf{R}$ and $C_2 \geq 0$ such that $\text{Ric}_{\omega_g} \geq - C_1 \omega_g + C_2 f^{\ast} \omega_h$, then \begin{eqnarray}\label{SchwarzInequality1}
\Delta_{\omega_g} | \partial f |^2 & \geq & -C_1 | \partial f |^2 + \left( \frac{\kappa_0}{r} + \frac{C_2}{n} \right) | \partial f |^4. 
\end{eqnarray} In particular, if $\kappa_0 n + r C_2 >0$,  we have \begin{eqnarray}\label{SchwarzLemma0}
| \partial f |^2 & \leq & \frac{C_1 r n}{\kappa_0 n + r C_2}.
\end{eqnarray} 

\hfill

\noindent The main application of \nameref{KahlerNewSL} is the following extension of the Wu--Yau theorem. 

\subsection*{Theorem 3.2}\label{MainWuYau}
Let $X$ be a compact K\"ahler manifold with a Hermitian metric satisfying $\text{RBC}_{\omega}^{\tau} \leq 0$ or a pluriclosed metric with $\text{HSC}_{\omega} \leq 0$.  Then the canonical bundle $K_X$ is nef.  If, in addition, $\text{RBC}_{\omega}^{\tau} <0$, or the metric is pluriclosed with $\text{HSC}_{\omega}<0$, then $X$ is projective with ample canonical bundle.  \begin{proof}
Let $X$ be a compact K\"ahler manifold with a Hermitian metric $\omega$ satisfying $\text{RBC}_{\omega}^{\tau} \leq 0$ or a pluriclosed metric $ \text{HSC}_{\omega} \leq 0$. We follow the argument in \cite{WuYau1, TosattiYang} and proceed by contradiction, supposing that $K_X$ is not nef. In particular,  $2\pi c_1(K_X^{-1})$ does not lie on the boundary of the K\"ahler cone.  Fix a K\"ahler metric $\omega_0$ on $X$, and let $\varepsilon_0>0$ be the constant such that $\varepsilon_0 \{ \omega_0 \} - 2\pi c_1(K_X^{-1})$ intersects the boundary of the K\"ahler cone.  Then for any $\varepsilon>0$, the cohomology class $(\varepsilon+\varepsilon_0) \{ \omega_0 \} - 2\pi c_1(K_X^{-1})$ is a K\"ahler class (i.e., is contained in the interior of the K\"ahler cone).  Hence, there is a smooth function $\varphi_{\varepsilon} \in \mathcal{C}^{\infty}(X, \mathbf{R})$ such that $(\varepsilon + \varepsilon_0) \omega_0 - \text{Ric}_{\omega_0} + \sqrt{-1} \partial \bar{\partial} \varphi_{\varepsilon}$ is a K\"ahler metric. By the Aubin--Yau theorem \cite{Aubin,Yau1976}, there is a function $\psi_{\varepsilon} \in \mathcal{C}^{\infty}(X, \mathbf{R})$ such that $$\omega_{\varepsilon} \ : = \ (\varepsilon + \varepsilon_0) \omega_0 - \text{Ric}_{\omega_0} + \sqrt{-1} \partial \bar{\partial} u_{\varepsilon}$$ solves the Monge--Amp\`ere equation $\omega_{\varepsilon}^n = e^{u_{\varepsilon}} \omega_0^n$,  where $u_{\varepsilon} : = \varphi_{\varepsilon} + \psi_{\varepsilon}$.  Since $X$ is compact, there is a uniform constant $C_0>0$ such that $C_0^{-1} \omega_0 \leq \omega \leq C_0 \omega_0$. Hence, by differentiating the complex Monge--Amp\`ere equation, we see that \begin{eqnarray*}
\text{Ric}_{\omega_{\varepsilon}} & = & - \omega_{\varepsilon} + (\varepsilon + \varepsilon_0) \omega_0 \ \geq \ - \omega_{\varepsilon} + C_0^{-1} (\varepsilon + \varepsilon_0) \omega. 
\end{eqnarray*}
Applying \nameref{KahlerNewSL} to the identity map $f= \text{id} : (X,\omega_{\varepsilon}) \to (X, \omega)$, with $C_1 = 1$, $C_2 = C_0^{-1}(\varepsilon + \varepsilon_0)$, and $\kappa_0=0$, we see that \begin{eqnarray*}
| \partial f |^2 \ = \ \text{tr}_{\omega_{\varepsilon}}(\omega) & \leq & \frac{n}{C_0^{-1}(\varepsilon+\varepsilon_0)}, 
\end{eqnarray*}
which is uniformly bounded above as $\varepsilon \searrow 0$. The proof of the higher-order estimates in \cite{WuYau1,TosattiYang} remains unchanged, and hence, by extracting a smoothly convergent subsequence, we obtain the desired contradiction. This shows that $K_X$ is nef.

Since $K_X$ is nef,  for any $\varepsilon>0$, there is a smooth function $u_{\varepsilon} \in \mathcal{C}^{\infty}(X, \mathbf{R})$ such that $$\omega_{u_{\varepsilon}} \ : = \ \varepsilon \omega_0 - \text{Ric}_{\omega_0} + \sqrt{-1} \partial \bar{\partial} u_{\varepsilon}$$ solves the complex Monge--Amp\`ere equation $\omega_{u_{\varepsilon}}^n = e^{u_{\varepsilon}} \omega_0^n$.  Hence, by differentiating the complex Monge--Amp\`ere equation, \begin{eqnarray*}
    \text{Ric}_{\omega_{u_{\varepsilon}}} &=& - \sqrt{-1} \partial \bar{\partial} u_{\varepsilon} + \text{Ric}_{\omega_0} \ = \ - \omega_{u_{\varepsilon}} + \varepsilon \omega_0 \ \geq \ - \omega_{u_{\varepsilon}} + \varepsilon C_0^{-1} \omega.
\end{eqnarray*} Applying \nameref{KahlerNewSL} to the identity map $f=\text{id} : (X, \omega_{u_{\varepsilon}}) \to (X, \omega)$ with $C_1 = 1$, $C_2 = \varepsilon C_0^{-1}$ and $\kappa_0 = \kappa_0 <0$ then yields $$| \partial f |^2 \ = \ \text{tr}_{\omega_{u_{\varepsilon}}}(\omega) \ \leq \ \frac{n}{C_0^{-1} \varepsilon + \kappa_0}, $$ which is uniformly bounded above as $\varepsilon \searrow 0$. As before, the higher-order estimates in  \cite{WuYau1,TosattiYang} then apply without change to obtain the desired contradiction, showing that $K_X$ is ample.
\end{proof}

\subsection*{Theorem 3.3}\label{Theorem1.7}
Let $(X, \omega)$ be a complete Hermitian manifold with $  \text{Ric}_{\omega}^{\tau} \geq - C \omega$,  bounded torsion, and $  \text{HBC}_{\omega} \leq - B <0$. Then $X$ is not biholomorphic to the polydisk $\mathbf{D}^n$ for $n>1$.

\begin{proof}
Proceed by contradiction and suppose that the polydisk $\mathbf{D}^n$ has a complete Hermitian metric with bounded torsion,  $ \text{Ric}_{\omega}^{\tau} \geq - A$ and $ \text{HBC}_{\omega} \leq - B <0$.  Following \cite{TangYang}, it suffices to consider the case $n=2$ since the argument readily extends without change to higher dimensions. Let $h$ be the Hermitian metric underlying $\omega$, and let $\rho$ denote the Poincar\'e metric on $\mathbf{D}$.  The Schwarz lemma \cite{Lu,YauSchwarzLemma} applied to the inclusion $\iota_z : (\mathbf{D},\rho) \to (\mathbf{D}^2, \omega_h)$, where $\iota_z(w) = (z,w)$,  implies that $\text{tr}_{\rho}(\iota_z^{\ast} \omega_h) \leq 4/B$. Hence,  $$h_{2\overline{2}}(z,w) \ \leq \ \frac{4}{B(1- |w |^2)^2}.$$ Let $f : \mathbf{D} \to \mathbf{R}$ be the function defined by $f(z) : = h_{2 \overline{2}}(z,0)$.  Then $f$ is a smooth, positive bounded function.  Computing its Laplacian, we see that \begin{eqnarray*}
\Delta_{\rho} f & = & (1- | z |^2)^2 \frac{\partial^2 h_{2\overline{2}}}{\partial  z \partial \overline{z}} \ = \ (1- | z |^2)^2 \left( - R_{1 \overline{1} 2 \overline{2}}(z,0) + h^{\alpha \overline{\beta}} \frac{\partial h_{2 \overline{\beta}}}{\partial z} \frac{\partial h_{\alpha \overline{2}}}{\partial  \overline{z}} \right).
\end{eqnarray*}
Since $h^{\alpha \overline{\beta}} \frac{\partial h_{2 \overline{\beta}}}{\partial z} \frac{\partial h_{\alpha \overline{2}}}{\partial \overline{z}} \geq 0$ and $R_{1 \overline{1} 2 \overline{2}}(z,0) \leq - B h_{1 \overline{1}}(z,0) h_{2 \overline{2}}(z,0)$, it follows that \begin{eqnarray*}
\Delta_{\rho} f & \geq & B(1- | z |^2)^2 h_{1 \overline{1}}(z,0) h_{2 \overline{2}}(z,0).
\end{eqnarray*} The Schwarz lemma applied to the projection $\pi : (\mathbf{D}^2, \omega_h) \to (\mathbf{D},\omega)$, where $\pi(z,w) = z$, yields $$(1- | z |^2)^2 h_{1 \overline{1}}(z,w) \ \geq \ \frac{4}{A},$$ and therefore, $\Delta_{\rho} f  \geq  \frac{4B}{A} f$. The desired contradiction is obtained from the Omori--Yau maximum principle \cite{Omori, YauMaximumPrinciple}.
\end{proof}

\section{The Tempered Hermitian Curvature Flows} \label{sec:Flows}
To prove a general Wu--Yau theorem void of any ambient K\"ahler assumption,  the only known method is to use a geometric flow that incorporates the second Chern Ricci curvature \cite{LeeHCF, LeeStreets}.  Streets--Tian \cite{StreetsTian} introduced a family of \textit{Hermitian curvature flows} given by $$\frac{\partial \omega_t}{\partial t} \ = \ -  \text{Ric}_{\omega_t}^{(2)} + \mathcal{Q}_{\omega_t},$$ where $\mathcal{Q}_{\omega_t}$ is an arbitrary quadratic expression in the torsion of the Chern connection. The flows in this family are strictly parabolic, and hence always admit unique short-time solutions on compact manifolds. 

The pluriclosed flow was used by Lee--Streets \cite{LeeStreets} to show that a compact pluriclosed manifold with a Hermitian metric of negative real bisectional curvature is projective with ample canonical bundle. The key observation is that a negative upper bound on the real bisectional curvature of a background metric gives a uniform $C^0$--bound for metrics evolving under the pluriclosed flow via a parabolic Schwarz lemma. It is natural to suspect that using the same method, the Lee--Streets theorem could be improved, requiring only a negative upper bound on the tempered real bisectional Curvature $ \text{RBC}_{\omega}^{\tau}<0$.  This does not appear to be possible,  since the quadratic torsion term that tempers the second Chern Ricci curvature does not coincide with the quadratic torsion term that appears in the pluriclosed flow. In light of this, we now give a parabolic version of the tempered Schwarz lemma:

\subsection*{Theorem 4.1}\label{ParabolicSchwarz}
Let $X$ be a complex manifold endowed with a smooth family of Hermitian metrics $(\omega_t)_{t \in I}$ satisfying
\begin{eqnarray}\label{SuperSoln}
\frac{\partial}{\partial t}\omega_t & \geq & -\text{Ric}^\tau_{\omega_t} - \omega_t,
\end{eqnarray} for some fixed $\tau \in \mathbf{R}_{>0}$.  
If $\omega_h$ is a metric on $X$ with $\text{RBC}_{\omega_h}^\tau \leq -\kappa_0 < 0$, then
 \begin{eqnarray*}
\left( \frac{\partial}{\partial t} - \Delta_{\omega_t} \right) \text{tr}_{\omega_t}(\omega_h) & \leq & - \frac{\kappa_0}{n} \text{tr}_{\omega_t}(\omega_h)^2 + \text{tr}_{\omega_t}(\omega_h),
\end{eqnarray*}
In particular, if $X$ is compact, then $\sup_{t \in I}\text{tr}_{\omega_t}(\omega_h) \leq \max\left\{\sup_X \text{tr}_{\omega_0}(\omega_h) , \frac{\kappa_0}{n}\right\}$.
\begin{proof}
This is a straightforward extension of \nameref{MainTheoremTemperedRBC}.  Here are the details: \begin{eqnarray*}
\frac{\partial}{\partial t} \text{tr}_{\omega_t}(\omega_h) \ = \ \frac{\partial}{\partial t} g^{k \overline{\ell}} h_{k \overline{\ell}} &=& - g^{p \overline{\ell}} g^{k \overline{q}} h_{k \overline{\ell}} \left( \frac{\partial }{\partial t} g_{p \overline{q}} \right) \\
&\leq& - g^{p \overline{\ell}} g^{k \overline{q}} h_{k \overline{\ell}} \left( -  \text{Ric}^{(2)}_{p \overline{q}} - \frac{1}{4}(1-\tau^{-1}) \mathcal{Q}_{p \overline{q}}^2 - g_{p \overline{q}} \right).
\end{eqnarray*}
From \nameref{MainTheoremTemperedRBC}, we have \begin{eqnarray*}
\Delta_{\omega_t} \text{tr}_{\omega_t}(\omega_h) & \geq &   g^{p \overline{\ell}} g^{k \overline{q}} h_{k \overline{\ell}}   \text{Ric}^{(2)}_{p \overline{q}} - \frac{1}{4}(1-\tau^{-1}) \mathcal{Q}_{p \overline{q}}^2 g^{p \overline{\ell}} g^{k \overline{q}} h_{k \overline{\ell}} \\
&& \hspace{3cm} -  R_{i \overline{j} k \overline{\ell}} g^{i \overline{j}} g^{k \overline{\ell}} + \frac{1}{4}(1-\tau) \mathcal{Q}_{i \overline{j} k \overline{\ell}} g^{i \overline{j}} g^{k \overline{\ell}}.
\end{eqnarray*}
Hence,  \begin{eqnarray*}
\left( \frac{\partial}{\partial t} - \Delta_{\omega_t} \right) \text{tr}_{\omega_t}(\omega_h) & \leq & -\frac{\kappa_0}{n} \text{tr}_{\omega_t}(\omega_h)^2 + \text{tr}_{\omega_t}(\omega_h),
\end{eqnarray*}
and the final claim follows from the maximum principle.
\end{proof}

As remarked previously, metrics evolving under (normalized) pluriclosed flow do not necessarily satisfy \eqref{SuperSoln}. It is therefore natural to consider the (normalized) \emph{tempered Hermitian curvature flows} defined as supersolutions of the following special case of the Hermitian curvature flow:\begin{eqnarray} \label{THCF}
\frac{\partial \omega_t}{\partial t} &=& -  \text{Ric}_{\omega_t}^{(2)} - \frac{1}{4}(1-\tau^{-1}) \mathcal{Q}_{\omega_t}^2 - \omega_t,
\end{eqnarray}
with $\omega_t \vert_{t=0} = \omega_0$.  By \nameref{ParabolicSchwarz}, solutions to \eqref{THCF} satisfy a parabolic Schwarz lemma. If an  understanding of the long-time existence and solutions of the tempered Hermitian curvature flows can be achieved, this may provide a significantly more general extension of the theorem of Lee--Streets \cite{LeeStreets}:

\subsection*{Question 4.2}
Let $(X, \omega)$ be a compact Hermitian manifold with $ \text{RBC}_{\omega}^{\tau}<0$.  Does the tempered Hermitian curvature flow exist for all time? Does it converge to a K\"ahler current?\\

Indeed, as in \cite{LeeStreets}, if the tempered Hermitian curvature flow converges to a K\"ahler current,  then \cite{DemaillyPaun} implies that $X$ is in the Fujiki class $\mathcal{C}$ (i.e., bimeromorphic to a compact K\"ahler manifold). Since $ \text{RBC}_{\omega}^{\tau}  <0$, every holomorphic map $\mathbf{P}^1 \to X$ is constant,  and hence, $X$ is K\"ahler.  Thus, by applying \nameref{MainWuYau}, we see that $X$ is projective with ample canonical bundle.

On the other hand, if \nameref{ParabolicSchwarz} can be adapted to apply to the pluriclosed flow, then \nameref{thm:main} can be extended to the following: A compact pluriclosed manifold $(X, \omega)$ with $\text{HSC}_{\omega} <0$ is projective with ample canonical bundle.

\section{The Schwarz Lemma for Gauduchon Connections} \label{sec:Gauduchon}
In this final section, we discuss the dependence of the Hermitian connection in the Schwarz lemma and extend the Schwarz lemma to the Gauduchon connections. The affine line \eqref{GaudFormula} of Gauduchon connections arises naturally from representation-theoretic considerations (see \cite{GauduchonHermitianConnections, BroderStanfield,StanfieldThesis}). For $t = -1$, the Gauduchon connection ${}^t \nabla$ recovers the \textit{Strominger--Bismut connection} ${}^b \nabla$, defined by its torsion being totally skew-symmetric. The Strominger--Bismut connection plays an important role in heterotic string theory and the pluriclosed flow \cite{StreetsTianPCF} (see, e.g., \cite{BroderStanfield} and the references therein). 

To understand the role played by the choice of connection in the Schwarz lemma estimate, and to establish the Schwarz lemma for the Gauduchon connections, we introduce the following ad-hoc definition.

\subsection*{Definition 5.1}\label{Def52}
Let $(X, \omega)$ be a Hermitian manifold and take $\nabla$ to be a Hermitian connection on $T^{1,0}X$. We say that $\nabla$ is \textit{suitably skew-symmetric} if the torsion ${}^{\nabla} T$ of $\nabla$ satisfies \begin{equation}\label{eqn:skew10}
	g\left({}^\nabla T(\overline u,v),\overline w\right) = - g\left({}^\nabla T(\overline w,v),\overline u\right),
\end{equation}
for all $u,v,w \in T^{1,0}X$. \\

The Gauduchon connections ${}^t \nabla$ are suitably skew-symmetric (c.f., \cite[Lemma 2.7]{BroderStanfield}). More generally, for any smooth function $u : X \to \mathbf{C}$, the Hermitian connection ${}^u \nabla$ defined by \begin{eqnarray*}
    g\left({}^{u} T(\overline u,v),\overline w\right) = \frac{1-u}{2}g\left(v,\overline{{}^cT(u,w)}\right),
\end{eqnarray*} is suitably skew-symmetric, where ${}^c T$ is the torsion of the Chern connection.

The following shows the precise effect of changing the Hermitian connection (amongst the suitably skew-symmetric Hermitian connections) in the Schwarz lemma:

\subsection*{Proof of \nameref{thm:OtherConnections}}
It suffices to show that 
\[
 \text{Sym}\left( {}^{\mathcal{T}} \nabla^{1,0} \partial f \right) \ = \ \text{Sym}\left( {}^{c} \nabla \partial f \right).
\]
The Hermitian connection $\nabla$ is determined by the $(1,1)$--part of its torsion (see, e.g., \cite[\S 2.3]{BroderStanfield}), given by
\[
g\left(\nabla_{u}v,\overline{w}\right) \ = \ g\left({}^c\nabla_{u}v,\overline w \right) - g\left(v,{}^\nabla T(u, \overline w)\right), 
\] where $u,v,w \in T^{1,0}X$.  Since $\nabla$ is suitably skew-symmetric, for $u,v \in T^{1,0}X$, we have 
\[
\nabla_uv + \nabla_vu \ = \ {}^c \nabla_{u}v + {}^c\nabla_{v}u.
\]
The same formula holds for $\widetilde \nabla$ with the obvious changes. Hence, for all $u,v\in T^{1,0}X$ and $w \in T^{1,0}Y$,
\begin{eqnarray*}
    2h\left(\text{Sym}\left({}^\mathcal{T}\nabla^{1,0} \partial f \right)(u,v),\overline w\right) & = & h
    \left({}^\mathcal{T}\nabla_u(\partial f)(v) + {}^\mathcal{T}\nabla_v(\partial f)(u),\overline w\right)\\
    & = & h\left(\widetilde\nabla_{\partial f (u)}\left(\partial f(v)\right) + \widetilde\nabla_{\partial f (v)}\left(\partial f(u)\right),\overline w\right) \\
    &&\quad - h\left(\partial f\left(\nabla_uv + \nabla_v u\right),\overline {w}\right)\\
    & = & h\left({}^c\nabla_{\partial f (u)}\left(\partial f(v)\right) + {}^c\nabla_{\partial f (v)}\left(\partial f(u)\right),\overline w\right) \\
    &&\quad - h\left(\partial f\left({}^c\nabla_uv + {}^c\nabla_v u\right),\overline {w}\right)\\
    &=& 2h\left(\text{Sym}\left({}^c\nabla^{1,0} \partial f\right),\overline w \right),
\end{eqnarray*}
and the claim follows since $f$ is holomorphic. \hfill $\Box$

\hfill

From \nameref{thm:OtherConnections}, we can obtain the Schwarz lemma for the Gauduchon connections, extending the known Schwarz lemma for the Chern connection. To this end, we require the following lemma.

\subsection*{Lemma 5.2}\label{lem:GaudRandT}
Let $(X, \omega_g)$ be a Hermitian manifold. Let ${}^t T$ and ${}^t R$ denote the torsion and curvature of the Gauduchon connection ${}^t \nabla$. Then, in any local unitary frame, we have 
\begin{itemize}
	\item[(i)] ${}^tT_{ij}^k = t {}^c T_{ij}^k$, for all $t \in \mathbf{R}$, and
	\item[(ii)] for $t \in \mathbf{R} \backslash \left \{ 0,\frac{1}{2} \right \}$,
 \begin{eqnarray*}
 {}^c R_{i\overline j k \overline \ell} &=& \frac{t^2 + 2t -1}{2t(2t-1)}{}^tR_{i\overline j k \overline \ell } + \frac{(t-1)^2}{2t(2t-1)}{}^tR_{k \overline \ell i \overline j} + \frac{t-1}{2(2t-1)}\left({}^tR_{k\overline j i \overline \ell} + {}^tR_{i \overline \ell k \overline j}\right)\\
&&- \frac{(t-1)^2}{4t^2(2t-1)}{}^t T^r_{ik}\overline{{}^tT^r_{j\ell}} + \frac{(t-1)^2(t^2+2t-1)}{8t^3(2t-1)}{}^t T^{\ell}_{ir}\overline{{}^tT^k_{jr}} + \frac{(t-1)^4}{8t^3(2t-1)}{}^tT^j_{kr}\overline{{}^tT^i_{\ell r}}\\
&& + \frac{(t-1)^3}{8t^2(2t-1)}\left({}^tT^{\ell}_{kr}\overline{{}^t T^i_{jr}} + {}^t T^j_{ir}\overline{{}^t T^k_{\ell r}}\right).
\end{eqnarray*}
\end{itemize}
\begin{proof}
Let $\{e_i\}_{i=1}^n$ be a unitary frame. Recall that ${}^\ell \nabla$ is the unique Hermitian connection whose torsion satisfies ${}^\ell T(e_i,e_j) = 0$ for all $1\leq i,j \leq n$ (see \cite{GauduchonHermitianConnections}). Hence, by Equation \eqref{GaudFormula}, ${}^tT_{ij}^k = t{}^cT_{ij}^k$, which proves (i). For (ii), we recall from \cite[Theorem 2.8]{BroderStanfield} that 
	\[
{}^tR_{i\overline j k \overline \ell} \ = \  t{}^cR_{i\overline j k \overline \ell } + \frac{1-t}{2}\left({}^cR_{k\overline j i \overline \ell } + {}^cR_{i\overline \ell k \overline j}\right) + \left(\frac{1-t}{2}\right)^2 \left({}^cT^r_{ik}\overline{{}^cT^r_{j \ell}} - {}^cT^{\ell}_{ir}\overline{{}^cT^k_{jr}}\right).
\]
Hence, 
\begin{eqnarray*}
{}^tR_{i\overline j k \overline \ell} + {}^tR_{k \overline \ell i\overline j} &=& t\left({}^cR_{i\overline j k \overline \ell} + {}^cR_{k \overline \ell i \overline j}\right) + (1-t)\left({}^cR_{k\overline j i \overline \ell} + {}^cR_{i\overline \ell k \overline j}\right) \\
&& \hspace{3cm} + \left(\frac{1-t}{2}\right)^2 \left(2{}^cT^r_{ik}\overline{{}^cT^r_{j\ell}} - {}^cT^{\ell}_{ir}\overline{{}^cT^k_{jr}} - {}^cT^j_{kr}\overline{{}^cT^i_{\ell r}}\right).
\end{eqnarray*}
Interchanging $i$ and $k$ yields
\begin{eqnarray*}
{}^tR_{k\overline j i \overline \ell } + {}^tR_{i \overline \ell  k\overline j} &=& t\left({}^cR_{k\overline j i \overline \ell } + {}^cR_{i \overline \ell  k \overline j}\right) + (1-t)\left({}^cR_{i\overline j k \overline \ell } + {}^cR_{k\overline \ell  i \overline j}\right) \\
&& \hspace{3cm} + \left(\frac{1-t}{2}\right)^2 \left(2{}^cT^r_{ki}\overline{{}^cT^r_{j\ell }} - {}^cT^{\ell}_{kr}\overline{{}^cT^i_{jr}} - {}^cT^j_{ir}\overline{{}^cT^k_{\ell r}}\right).
\end{eqnarray*}
These furnish the linear system
\[
\begin{pmatrix}
t&1-t\\
1-t&t
\end{pmatrix}\begin{pmatrix}
{}^cR_{i\overline j k \overline \ell } + {}^cR_{k \overline \ell i \overline j}\\
{}^cR_{k\overline j i \overline \ell } + {}^cR_{i\overline \ell k \overline j}
\end{pmatrix} = \begin{pmatrix}
{}^tR_{i\overline j k \overline \ell } + {}^tR_{k \overline \ell  i \overline j} - \left(\frac{1-t}{2}\right)^2 \left(2{}^cT^r_{ik}\overline{{}^cT^r_{j\ell }} - {}^cT^{\ell}_{ir}\overline{{}^cT^k_{jr}} - {}^cT^j_{kr}\overline{{}^cT^i_{\ell r}}\right)\\
{}^tR_{k\overline j i \overline \ell } + {}^tR_{i \overline \ell  k \overline j} - \left(\frac{1-t}{2}\right)^2 \left(-2{}^cT^r_{ik}\overline{{}^cT^r_{j\ell }} - {}^cT^{\ell}_{kr}\overline{{}^cT^i_{jr}} - {}^cT^j_{ir}\overline{{}^cT^k_{\ell r}}\right)
\end{pmatrix},
\]
and hence,
\begin{eqnarray*}
{}^cR_{k\overline j i \overline \ell } + {}^cR_{i\overline \ell  k \overline j} &=& \frac{t-1}{2t-1}\left({}^tR_{i\overline j k \overline \ell } + {}^tR_{k \overline \ell i \overline j}\right) + \frac{t}{2t-1}\left({}^tR_{k\overline j i \overline \ell } + {}^tR_{i \overline \ell k \overline j}\right)\\
&&- \frac{(t-1)^3}{4(2t-1)}\left(2{}^cT^r_{ik}\overline{{}^cT^r_{j\ell}} - {}^cT^{\ell}_{ir}\overline{{}^cT^k_{jr}} - {}^cT^j_{kr}\overline{{}^cT^i_{\ell r}}\right) \\
&&- \frac{t(t-1)^2}{4(2t-1)}\left(-2{}^cT^r_{ik}\overline{{}^cT^r_{j\ell}} - {}^cT^{\ell}_{kr}\overline{{}^cT^i_{jr}} - {}^cT^j_{ir}\overline{{}^cT^k_{\ell r}}\right).
\end{eqnarray*}

\noindent As a consequence, we have
\begin{eqnarray*}
t{}^cR_{i\overline j k \overline \ell } &= &  {}^tR_{i\overline j k \overline \ell } - \frac{(t-1)^2}{4}\left({}^cT^r_{ik}\overline{{}^cT^r_{j\ell }} - {}^cT^{\ell}_{ir}\overline{{}^cT^k_{jr}}\right)\\
&& \hspace{1cm} +\frac{(t-1)^2}{2(2t-1)}\left({}^tR_{i\overline j k \overline \ell } + {}^tR_{k \overline \ell  i \overline j}\right) + \frac{t(t-1)}{2(2t-1)}\left({}^tR_{k\overline j i \overline \ell } + {}^tR_{i \overline \ell k \overline j}\right)\\
&& \hspace{2cm} - \frac{(t-1)^4}{8(2t-1)}\left(2{}^cT^r_{ik}\overline{{}^cT^r_{j\ell}} - {}^cT^{\ell}_{ir}\overline{{}^cT^k_{jr}} - {}^cT^j_{kr}\overline{{}^cT^i_{\ell r}}\right) \\
&& \hspace{3cm} - \frac{t(t-1)^3}{8(2t-1)}\left(-2{}^cT^r_{ik}\overline{{}^cT^r_{j\ell}} - {}^cT^{\ell}_{kr}\overline{{}^cT^i_{jr}} - {}^cT^j_{ir}\overline{{}^cT^k_{\ell r}}\right)\\
&& = \frac{t^2 + 2t -1}{2(2t-1)}{}^tR_{i\overline j k \overline \ell } + \frac{(t-1)^2}{2(2t-1)}{}^tR_{k \overline \ell i \overline j} + \frac{t(t-1)}{2(2t-1)}\left({}^tR_{k\overline j i \overline \ell } + {}^tR_{i \overline \ell k \overline j}\right)\\
&& \hspace{1cm} - \frac{t(t-1)^2}{4(2t-1)}{}^cT^r_{ik}\overline{{}^cT^r_{j\ell}} + \frac{(t-1)^2(t^2+2t-1)}{8(2t-1)}{}^cT^{\ell}_{ir}\overline{{}^cT^k_{jr}} \\
&& \hspace{2.5cm} + \frac{(t-1)^4}{8(2t-1)}{}^cT^j_{kr}\overline{{}^cT^i_{\ell r}} + \frac{t(t-1)^3}{8(2t-1)}\left({}^cT^{\ell}_{kr}\overline{{}^cT^i_{jr}} + {}^cT^j_{ir}\overline{{}^cT^k_{\ell r}}\right).
\end{eqnarray*}
The result now follows from item (i).
\end{proof}

Recall that in \cite{BroderTangAltered,BroderStanfield} we defined the \emph{$t$--Gauduchon altered real bisectional curvature} by
$${}^t\widetilde{\text{RBC}}_{\omega}(\xi) \ := \ \frac{1}{| \xi |_{\omega}^2}  \sum_{\alpha, \beta, \gamma, \delta} {}^t R_{\alpha \overline \delta \gamma \overline \beta} \xi^{\alpha \overline \beta}\xi^{\gamma \overline \delta},$$ for all non-negative Hermitian $(1,1)$--tensors $\xi$.  From \nameref{lem:GaudRandT} we inherit formulae for the tempered Ricci curvature ${}^c \text{Ric}^{\tau}$ and tempered real bisectional curvature ${}^c \text{RBC}^{\tau}$ in terms of the curvature and torsion of the Gauduchon connections:

\subsection*{Proposition 5.3}\label{Prop44}
For all $t\in \mathbf{R} \backslash \left \{ 0, \frac{1}{2} \right \}$, and any local unitary frame $e^k$, we have 
\begin{eqnarray*}
{}^c \text{Ric}_{k \overline{\ell}}^{\tau} &=& \frac{t^2 + 2t -1}{2t(2t-1)}{}^t\text{Ric}_{k \overline{\ell}}^{(2)} + \frac{(t-1)^2}{2t(2t-1)}{}^t\text{Ric}_{k \overline{\ell}}^{(1)} + \frac{t-1}{2(2t-1)}\left({}^t\text{Ric}_{k \overline{\ell}}^{(3)} + {}^t\text{Ric}_{k \overline{\ell}}^{(4)}\right)\\
&& + \frac{(t-1)^2\left(t^2 - 4t + 1\right)}{8t^3(2t-1)} {}^t T^r_{ik}\overline{{}^tT^r_{i\ell }}+ \frac{(t-1)^3}{4t^2(2t-1)}\text{Re}{}^tT^{\ell}_{kr}\overline{{}^tT^i_{ir}}\\
&& +\left(\frac{(t-1)^2(t^2+2t-1)}{8t^3(2t-1)} + \frac{1-\tau^{-1}}{4t^2}\right){}^tT^{\ell}_{ir}\overline{{}^tT^k_{ir}}
\end{eqnarray*}
and
\begin{eqnarray*}
{}^c \text{RBC}^\tau(\xi) &=& \frac{t}{2t-1}  \phantom{.} {}^t\text{RBC}(\xi) + \frac{(t-1)}{2t-1}  {}^t\widetilde{\text{RBC}}(\xi) -\left(\frac{(t-1)^2}{4t^2(2t-1)} + \frac{1-\tau}{4t^2}\right){}^t T_{ik}^r\overline{{}^tT_{j\ell}^r} \xi^{i\overline j}\xi^{k \overline \ell} \\
&& +  \frac{t(t-1)^2}{4t^2(2t-1)}{}^tT^{\ell}_{ir}\overline{{}^t T^k_{jr}} \xi^{i\overline j}\xi^{k \overline \ell}  +\frac{(t-1)^3}{4t^2(2t-1)}{}^t T^j_{ir}\overline{{}^tT^k_{\ell r}} \xi^{i\overline j} \xi^{k \overline \ell},
\end{eqnarray*}
for all nonnegative Hermitian $(1,1)$--tensors $\xi$. \\

Recall that the Strominger--Bismut connection ${}^b \nabla$ corresponds to the Gauduchon parameter $t=-1$. By evaluating each of the rational functions in \nameref{Prop44} at $t=-1$, we have the following Schwarz lemma for the Strominger--Bismut connection:

\subsection*{Corollary 5.4}
Let $f : (X, \omega_g) \to (Y, \omega_h)$ be a holomorphic map of rank $r$ between Hermitian manifolds. Endow $T^{1,0}X$ and $T^{1,0}Y$ with the Strominger--Bismut connection ${}^b \nabla$ and ${}^b \widetilde{\nabla}$, respectively. Then \begin{eqnarray*}
\Delta_{\omega_g} | \partial f |^2 & \geq &  \frac{1}{3} \left( 2 {}^b \text{Ric}_{k \overline{\ell}}^{(1)} -{}^b \text{Ric}_{k \overline{\ell}}^{(2)}  +   {}^b \text{Ric}_{k \overline{\ell}}^{(3)} + {}^b \text{Ric}_{k \overline{\ell}}^{(4)} \right) g^{k \overline{q}} g^{p \overline{\ell}} f_p^{\alpha} \overline{f_q^{\beta}} h_{\alpha \overline{\beta}} \\
&& \hspace{1cm} + \left( {}^b T^r_{ik}\overline{{}^b T^r_{i\ell}}+ \frac{2}{3} \text{Re}{}^b T^{\ell}_{kr}\overline{{}^b T^i_{ir}}  -\left( \frac{3+\tau}{12\tau}  \right){}^b T^{\ell}_{ir}\overline{{}^b T^k_{ir}} \right) g^{k \overline{q}} g^{p \overline{\ell}} f_p^{\alpha} \overline{f_q^{\beta}} h_{\alpha \overline{\beta}}  \\
&& \hspace{3cm} + \frac{1}{3} \left( {}^b R^h_{\alpha \overline{\beta}\gamma\overline{\delta}} + 2 {}^b R^h_{\alpha \overline{\delta} \gamma \overline{\beta}} \right) g^{i \overline{j}} f_i^{\alpha} \overline{f_j^{\beta}} g^{p \overline{q}} f_p^{\gamma} \overline{f_q^{\delta}} \\
&& + \frac{1}{3} \left( {}^b \tilde{T}^{\delta}_{\alpha \mu}\overline{{}^b \tilde{T}^{\gamma}_{\beta \mu}} + 2{}^b \tilde{T}^{\beta}_{\alpha \mu}\overline{{}^b \tilde{T}^{\gamma}_{\delta \mu}} + \left( 1 - \frac{(1-\tau)}{12} \right) {}^b \tilde{T}_{\alpha \gamma}^{\mu}\overline{{}^b \tilde{T}_{\beta \delta }^{\mu}}   \right) g^{i \overline{j}} f_i^{\alpha} \overline{f_j^{\beta}} g^{p \overline{q}} f_p^{\gamma} \overline{f_q^{\delta}}.
\end{eqnarray*}

\end{document}